\documentclass[12pt]{amsart}
\usepackage{amssymb,amsmath,latexsym,amscd,amsfonts,amsrefs}
\usepackage{graphics}
\usepackage{epsfig}
\usepackage{psfrag}
\usepackage{enumerate}
\def\ignore #1 {}

\newtheorem{thm}{Theorem}

\newtheorem{lem}[thm]{Lemma}
\newtheorem{cor}[thm]{Corollary}

\theoremstyle{definition}
\newtheorem{dfn}[thm]{Definition}
\newtheorem{rem}[thm]{Remark}
\newtheorem{ex}[thm]{Example}

\def\hpic #1 #2 {\mbox{$\begin{array}[c]{l} \epsfig{file=#1,height=#2} \end{arr\
ay}$}}
\def\vpic #1 #2 {\mbox{$\begin{array}[c]{l} \epsfig{file=#1,width=#2} \end{arra\
y}$}}

\def\N{{\mathbb N}}

\def\R{{\mathbb R}}

\def\G{{\mathbb G}}
\def\A{{\mathbb A}}

\def\S{{\mathbb S}}
\def\uS{{\underline{\mathbb S}}}
\def\H{{\mathbb H}}

\def\cG{{\mathcal G}}
\def\e{{\bf e}}

\begin{document}

\title{Arithmetical Semirings}

\author{Marco Aldi}

\address{Department of Mathematics and Applied Mathematics, Virginia Commonwealth University, United States.}

\begin{abstract}
We study the number of connected graphs with $n$ vertices that cannot be written as the cartesian product of two graphs with fewer vertices. We give an upper bound which implies that for large $n$ almost all graphs are both connected and cartesian prime. For graphs with an even number of vertices, a full asymptotic expansion is obtained. Our method, inspired by Knopfmacher's theory of arithmetical semigroups, is based on reduction to Wright's asymptotic expansion for the number of connected graphs with $n$ vertices.
\end{abstract}


\maketitle

\section{introduction}

Let $\mathbb G(n)$ be the number of all unlabeled finite simple graphs (here simply referred to as {\it graphs}) with $n$ vertices and let $\mathbb G^+(n)$ be the number of connected graphs with $n$ vertices. It is a classical result of Wright \cite{W2} that almost all graphs are connected or, more precisely, 
\begin{equation}\label{eq:-1}
\lim_{n\to +\infty} \frac{\G^+(n)}{\G(n)} = 1\,.
\end{equation}
Let $\G^\Box(n)$ be the number of connected graphs with $n$ vertices that cannot be written as the cartesian product of connected graphs with fewer vertices. In this paper we give a short proof that for $n$ sufficiently large
\begin{equation}\label{eq:0}
\mathbb G^+(n)-\mathbb G^\Box(n) \le \G\left(\left\lfloor\frac{n}{2}\right\rfloor \right)+\G\left(\left\lfloor\frac{n}{3}\right\rfloor+3 \right)\,, 
\end{equation}
which easily implies 
\begin{equation}\label{eq:000}
\lim_{n\to +\infty} \frac{\G^\Box(n)}{\G^+(n)} = 1\,.
\end{equation}
Informally, \eqref{eq:-1} and \eqref{eq:000} can be summarized as saying that almost all graphs are both connected and cartesian primes, a statement which appears to be well-known to experts, but that we were unable to locate in the literature. Furthermore, for graphs with an even number of vertices we obtain the full asymptotic expansion
\begin{equation}\label{eq:00}
\G^+(2n)-\G^\Box(2n) = \frac{2^{\binom{n}{2}}}{n!}+\sum_{s=1}^{R-1}\omega_s(n)\frac{2^{\binom{n-s}{2}}}{(n-s)!} + O\left(\frac{2^{\binom{n-R}{2}}}{(n-2R)!}\right)
\end{equation}
where $R$ is an arbitrary integer greater than $1$ and, for each integer $s\ge 1$, $\omega_s$ is an explicitly computable polynomial of degree $s$. 

Enumerative formulas for the number of cartesian prime graphs were derived in \cite{jili} in the context of the theory of species. In this paper we derive asymptotic information in the context of Knopfmacher's abstract analytic number theory \cite{K2}. The connection between graph theory and abstract analytic number theory was first observed in \cite{K} where it was shown that Wright's  results on the asymptotic distribution of connected graphs \cite{W1} and \cite{W2} can be generalized to the study of more general additive arithmetical semigroups. In this paper we take the analysis of \cite{K} further by introducing the notion of {\it arithmetical semiring} which can be thought of an abstraction of the semiring of graphs in which addition and multiplication are given by the disjoint union and the cartesian product, respectively. The main feature of arithmetical semirings is that the compatibility of the two operations allows us to answer questions about the asymptotic distribution of multiplicative primes by reducing them to analogous questions about additive primes. In our main example of graphs, this allows to use the powerful results of \cite{W1} and \cite{W2} to prove \eqref{eq:-1}, \eqref{eq:0}, \eqref{eq:000}, and \eqref{eq:00}. To illustrate the potential fruitfulness of working with abstract arithmetical semirings we discuss the example of graphs with an even number of edges, which is nontrivial due to the presence of additional cartesian-prime graphs. In \cite{K2} Knopfmacher showed that additional results about the distribution of connected graphs can be obtained by looking at certain functions such as divisor counting functions and Euler-type functions that are natural from the point of view of abstract analytic number theory. In this paper we derive multiplicative analogues of these results by introducing the notion of submultiplicative function on arithmetical semirings.

\bigskip\noindent \textbf{Acknowledgments:} The author would like to thank Richard Hammack for patiently answering his questions on the cartesian product of graphs and the anonymous referees for helpful feedback.

\section{Arithmetical Semirings}

In this section we define the main algebraic notion of this paper, arithmetical semirings, by equipping additive arithmetical semigroups with a second, compatible, operation. We refer the reader to \cite{K2} for a systematic treatment of arithmetical semigroups, including numerous examples and applications to several areas of mathematics. Here we limit ourselves to the definitions needed to make the present paper self-contained (assuming only standard algebraic terminology, see e.g.\cite{Lang}). 

\begin{dfn}
\normalfont A {\it commutative monoid} is a triple $(\mathbb U,\bullet,\e_\bullet)$ consisting of a set $\mathbb U$ endowed with a commutative associative operation $\bullet$ with identity $\e_\bullet\in \mathbb U$. A {\it unique factorization monoid} is a monoid $(\mathbb U, \bullet, \e_\bullet)$ admitting a set of generators $\mathbb U^\bullet\subseteq \mathbb U\setminus \{\mathbf e_\bullet\}$ with respect to which every $G\in \mathbb U\setminus \{\mathbf e_\bullet\}$ can be uniquely (up to permutation of the factors) written as
\[
G=G_{i_1}\bullet G_{i_2}\bullet\ldots \bullet G_{i_n}
\]
for some (not necessarily distinct) $G_{i_1},\ldots,G_{i_n}\in \mathbb U^{\bullet}$. If this is the case, the elements of $\mathbb U^\bullet$ are referred to as $\bullet$-primes. 
\end{dfn}

\begin{dfn}\label{dfn:additiveartihmeticasemigroup} 
\normalfont
Consider a commutative monoid $(\A,+,\mathbf e_+)$ that is a unique factorization monoid with respect to a set $\mathbb A^+$ of $+$-primes. Given a monoid homomorphism $\partial$ from $(\mathbb A,+,\e_+)$ to $(\R_{\ge 0},+,0)$, let $\underline \A$ and $\underline \A^+$ be the set-valued functions on $\mathbb R_{\ge 0}$ such that  $\underline{\A}(x)=\partial^{-1}(x)$ and $\underline{\A}^+(x)=\partial^{-1}(x)\cap \A^+$ for each $x\in \mathbb R_{\ge 0}$. We also define integer-valued functions $\A$ and $\A^+$ on $\mathbb R_{\ge 0}$ such that $\A(x)=|\underline{\mathbb A}(x)|$ and  $\A^+(x)=|\underline\A^+(x)|$ for each $x\in \R_{\ge 0}$. We say that $(\A,+,\partial)$ is an {\it additive arithmetical semigroup} if 
\begin{enumerate}[i)]
\item $\underline{\A}(0)=\{\bf e_+\}$;
\item $\A(x)=0$ unless $x$ is an integer;
\item $\A(x)$ is finite for each $x\in \R_{>0}$.

\end{enumerate}
If this is the case the homomorphism $\partial$ is called the {\it degree map} of the additive arithmetical semigroup. 
\end{dfn}

\begin{rem}\label{rem:3}
The reader may wonder why the additive monoid of reals is used in Definition \ref{dfn:additiveartihmeticasemigroup}, even though condition ii) implies that the degree map is integer valued. The reason is that some of the formulas in this paper are slightly more natural if, for suitable choices of additive arithmetical semigroup $\mathbb A$, one is allowed to write $0$ as $\mathbb A(x)$ for non-integer values of $x$.
\end{rem}

\begin{ex}\label{ex:graph-semigroup} Our main example is the additive arithmetical semigroup of graphs constructed as follows \cite{K}. Let $\G$ be the collection of all (unlabeled finite simple) graphs thought of as a semigroup with respect to the operation $+$ of disjoint union. It is straightforward to check that if the degree map $\partial$ is the function that to each graph assigns the corresponding number of vertices, then $(\G,+,\partial)$ is an additive arithmetical semigroup whose identity element is the empty graph. By definition, the set of additive primes $\G^+$ coincides with the set of {\it connected graphs}. 
\end{ex}

\begin{dfn}
A {\it commutative unital semiring} is a set $\S$ together with two operations $+$, $\Box$ such that
\begin{enumerate}[i)]
\item $(\S,+,\e_+)$ is a commutative monoid for some unit $\mathbf e_+\in \S$;
\item $(\S\setminus\{\mathbf e_+\},\Box, \e_\Box)$ is a commutative monoid for some unit $\mathbf e_\Box\in \S\setminus \{\mathbf e_+\}$;
\item $\Box$ is distributive with respect to $+$.
\end{enumerate} 
Let $(\S_1,+_1,\Box_1)$ and $(\S_2,+_2,\Box_2)$ be commutative unital semirings. A function $f:\S_1\to \S_2$ is a {\it semiring homomorphism} if $f(S+_1S')=f(S)+_2f(S')$ and $f(S\Box_1S')=f(S)\Box_2f(S')$ for all $S,S'\in \S$. 
\end{dfn}

\begin{dfn}\label{def:6}
Let $(\S,+,\Box)$ be a commutative unital semiring and let $\partial: (\S,+,\Box)\to (\R_{\ge 0}, + ,\cdot)$ be a semiring homomorphism. We say that $(\S,+,\Box,\partial)$ is an {\it arithmetical semiring} if $(\S,+,\partial)$ is an additive arithmetical semigroup and $(\S^+,\Box,\e_\Box)$ is a unique factorization monoid such that $\uS^{+}(1)=\{{\bf e}_\Box\}$. We refer to elements of $\S^+$ as {\it additive primes}. We write $\mathbb S^\Box$ for the set of {\it multiplicative primes} i.e.\ $\Box$-primes in the unique factorization monoid $(\S^+,\Box,\e_\Box)$. In analogy with Definition \ref{dfn:additiveartihmeticasemigroup}, for every real $x\ge 0$, we denote by $\uS^\Box(x)$ the set $\partial^{-1}(x)\cap \S^\Box$ of multiplicative primes of degree $x$ and write $\S^\Box(x)=|\uS^\Box(x)|$. 
\end{dfn}

\begin{rem}
If $(\S,+,\Box,\partial)$ is an arithmetical semiring, then in particular $(\S^+,\Box,\partial)$ is an arithmetical semigroup in the sense of \cite{K2}. We emphasize that Definition \ref{def:6} requires that the set of additive primes is closed under multiplication and that multiplicative primes form be a subset of the set of additive primes. These assumptions are used throughout the paper as a simple way to connect the asymptotic distribution of additive and multiplicative primes. In principle, useful asymptotic growth results might be obtained even if the multiplicative monoid does not have unique factorization or if the additive and multiplicative monoids do not necessarily share the same degree map. However, lacking of motivating examples, we choose to postpone the study of these more general structures to future work.
\end{rem}

\begin{ex}\label{ex:graph-semiring}
We now describe the arithmetical semiring or graphs which is our main example and our motivation for using the symbol $\Box$ to denote multiplication. Let $(\G,+,\partial)$ be the additive arithmetical semigroup of graphs as in Example \ref{ex:graph-semigroup}. The cartesian product (we refer the reader to \cite{HIK} for a systematic treatment of this notion) of two graphs $G_1,G_2\in \mathbb G$ is defined as the graph $G_1\Box G_2\in \mathbb G$ with vertices $V(G_1\Box G_2)=V(G_1)\times V(G_2)$ and edges
\begin{equation}\label{eq:cartesian}
E(G_1\Box G_2)=\tau(E(G_1)\times \Delta_2) \cup \tau(\Delta_1\times E(G_2))\,.
\end{equation}
Here $\Delta_i\subseteq (V(G_i))^2$ is the diagonal and $\tau:(V(G_1))^2\times (V(G_2))^2\to (V(G_1\Box G_2))^2$ is the involution interchanging the second and third components. According to a theorem of Sabidussi and Vizing (\cite{S}, \cite{V}), the commutative monoid of connected graphs $(\mathbb G^+,\Box)$ is a unique factorization monoid with multiplicative identity $\mathbf e_\Box=K_1$, the graph with only one vertex. Moreover, it is easy to check (see e.g.\ \cite{HIK}) that $(\G,+,\Box)$ is a unital semiring and that $\partial$ is a semiring homomorphism. Hence $(\G,+,\Box,\partial)$ is an arithmetical semiring. Multiplicative primes in $(\G,+,\Box,\partial)$ are known as {\it cartesian-prime graphs}. 
\end{ex}

\begin{rem}\label{rem:subsemiring}
It follows immediately from the definition that if $(\S,+,\Box,\partial)$ is an arithmetical semiring and $\S'\subseteq \S$ is a subsemiring that contains $\e_+$ and $\e_\Box$, then $(\S',+,\Box,\partial)$ is also an arithmetical semiring. 
\end{rem}

\begin{ex}\label{ex:subsemiring}
For instance it follows from \eqref{eq:cartesian} that the subset $\G_{\rm ev}\subseteq \G$ of graphs with an even number of edges is closed under the semiring operations defined in Example \ref{ex:graph-semiring}. Therefore, $(\G_{\rm ev},+,\Box,\partial)$ is an arithmetical semiring. 
\end{ex}

\begin{ex}\label{ex:hamming}
Another interesting subsemiring of $\G$ is the {\it Hamming semiring}
\[
\mathbb H=\N[K_2,K_3,\ldots]
\]
generated by all complete graphs (with constants counting connected components isomorphic to $K_1$). Therefore, $(\H,+,\Box,\partial)$ is an arithmetical semiring. The Hamming semiring allows us to introduce an infinite family of non-standard structures of arithmetical semiring on the semiring of graphs introduced in Example \ref{ex:graph-semiring}, as follows. Upon canonically identifying the Hamming semiring $\H$ with a semiring of polynomials in infinitely many variables over $\N$, the {\it Hamming polynomial} defined in \cite{hamming} can be thought of as the semiring homomorphism $c:(\G,+\Box)\to (\H,+,\Box)$ that to each graph $G\in \G$ assigns the graph-valued polynomial   
\[
c(G)=\sum \alpha_{i_1,\ldots,i_m} K_{i_1}\Box\cdots \Box K_{i_m}\,
\]
where $\alpha_{i_1,\ldots,i_m}$ (the indices $i_1,\ldots,i_m$ are possibly repeated) is the number of induced subgraphs of $G$ that are isomorphic to $K_{i_1}\Box\cdots \Box K_{i_m}$. Since $\partial(G)\le \partial(c(G))$, we conclude that $(\G,+,\Box,\partial\circ c)$ is an arithmetical semiring with underlying semiring $(\G,+,\Box)$. 
\end{ex}

\section{Asymptotic distribution of primes}

In this section we study the asymptotic distribution of additive and multiplicative primes in arithmetical semiring $(\S,+,\Box,\partial)$ whose underlying additive arithmetical semigroup satisfies certain growth conditions known as Knopfmacher's axioms \cite{K}. Our main observation is that the semiring structure forces the underlying multiplicative arithmetical semigroup to satisfy analogous growth conditions. In the case of graphs, we are able to extract precise asymptotic information about the distribution of cartesian-prime graphs from classical results about connected graphs.

\begin{rem}\label{rem:6}
Let $(\S,+,\Box)$ be an arithmetical semiring. By unique factorization, the assignment $S\mapsto S+\e_\Box$ defines an injection $\uS(n)\to \uS(n+1)$ for all $n$. In particular, the sequence $\{\S(n)\}$ is non-decreasing. However, the sequence $\{\S^+(n)\}$ that enumerates additive primes does not necessarily have the same property (see  Example \ref{ex:hamming2} below). For the purpose of comparing the asymptotic distribution of additive and multiplicative it is convenient to focus on arithmetical semirings that enjoy the following property.
\end{rem}

\begin{dfn}
An arithmetical semiring $(\S,+,\Box,\partial)$ is {\it monotonic} (resp.\ {\it strictly monotonic}) if there exists $N\in \N$ such that $\S^+(n)\le \S^+(n+1)$ (resp.\  $\S^+(n)< \S^+(n+1)$) for all $n>N$.
\end{dfn}

\begin{ex}
Let $(\G,+,\Box,\partial)$ be the arithmetical semiring of graphs defined in Example \ref{ex:graph-semiring}. Since for each graph $G\in \underline{\G}^+(n)$, the only vertex of degree $0$ in $G+K_1$ can be connected in at least one way to the vertices of $G$, then $(\G,+,\Box,\partial)$ is monotonic.
\end{ex}

\begin{ex}\label{ex:hamming2}
Let $(\H,+,\Box,\partial)$ be the Hamming arithmetical semiring introduced in Example \ref{ex:hamming}. Let $p$ be any prime number. By construction, the complete graph $K_p$ is the only connected graph with $p$ vertices in $\H$.  On the other hand, $K_{2p}$ and $K_2\Box K_p$ are both connected graphs in $\H$ with $2p$ vertices. Therefore, $1=\H^+(p')<\H^+(2p)$ for any two prime numbers $p,p'$. Since for every prime number $p$ there exists a prime number $p'$ such that $2p < p'$, the arithmetical semiring $(\H,+,\Box,\partial)$ is not monotonic. 
\end{ex}

\begin{rem}\label{rem:asymptoticnotation}
For the reminder of this paper we make systematic use of the following standard asymptotic notation. Given non-negative sequences $\{a_n\}$, $\{b_n\}$ we write $a_n=O(b_n)$ if there exist integers $N$, $C$ such that $a_n\le Cb_n$ for all $n\ge N$. Moreover, we write $a_n=o(b_n)$ if for every positive real number $\epsilon$ there exists an integer $N$ such that $a_n\le \epsilon b_n$ for all $n\ge N$. Finally, we write $a_n\sim b_n$ if $|a_n-b_n|=o(a_n)$.
\end{rem}

\begin{dfn}
Let $(\S,+,\Box,\partial)$ be an arithmetical semiring. We say  that 
\begin{enumerate}[1)]
\item {\it axiom $\cG_1^+$} holds if $\S(n)\sim \S^+(n)$; 
\item {\it axiom $\cG_1^\Box$} holds if $\S^+(n)\sim \S^\Box(n)$;
\item {\it axiom $\cG_2^+$} holds if $\S^+(n-1)=o(\S^+(n))$;
\item {\it axiom $\cG_2^\Box$} holds if  $\S^\Box\left(\left\lfloor \frac{n}{2}\right\rfloor\right)=o(\S^\Box(n))$;
\item {\it axiom $\cG_3^+$} holds if axiom $\cG_2^+$ holds and $\S(n)-\S^+(n)=O(\S^+(n-1))$ ;
\item {\it axiom $\cG_3^\Box$} holds if axiom $\cG_2^\Box$ holds and $\S^+(n)-\S^\Box(n)=O\left(\S^\Box\left(\left\lfloor \frac{n}{2}\right\rfloor\right) \right)$.
\end{enumerate}
\end{dfn}

\begin{rem}\label{rem:12}
The axioms $\mathcal G_i^+$, $k=1,2,3$ depend only on the underlying structure of additive arithmetical semigroup and are equivalent to the three axioms (in the respective order) defined in \cite{K}. The axioms $\cG_i^\Box$ are introduced here as their natural multiplicative analogues. 
\end{rem}

\begin{rem}
Roughly speaking, axiom $\cG_1^+$ states that almost all elements of the semiring $\S$ are additive primes. Similarly, if axiom $\cG_1^\Box$ holds, then in a sense almost all additive primes in $\S$ are also multiplicative primes. We emphasize that a choice of degree map is necessary in order to assign a rigorous meaning to the expression ``almost all'' and the validity of these axioms is a priori dependent on such a choice. Axioms $\cG_3^+$ and $\cG_3^\Box$ can be understood as more precise  versions of, respectively, axioms $\cG_1^+$ and $\cG_1^\Box$ in the sense that the comparison between the sequences involved in those axioms would become quantitative for a specific choice of the constants (denoted as $C$ and $N$ in Remark \ref{rem:asymptoticnotation}) that are implicit in the asymptotic notation. 
\end{rem}

\begin{rem}\label{rem:monotonic}
Axioms $\cG_2^+$ and $\cG_2^\Box$ can be thought of as a measure of how fast the number of additive and, respectively, multiplicative primes grows when compared to their degree. In particular, if $\cG_2^+$ holds the arithmetical semiring is strictly monotonic. The converse is not necessarily true and a priori strict monotonicity only implies $n=O(\S^+(n))$.
\end{rem}

\begin{lem}\label{lem:log}
Let $(\S,+,\Box,\partial)$ be an arithmetical semiring and suppose that there exist $\alpha,\beta,\gamma,a,b\in \R$ such that $\alpha>0$, $a>b>0$ and
\begin{equation}\label{eq:log}
\log \S(n) = \alpha n^{a+1} +\beta n \log n + \gamma n + O(n^b)\,.
\end{equation}
Then axioms $\cG_1^+$ and $\cG_2^+$ are satisfied. Moreover, if $B$ is the integer sequence defined recursively for all positive integers by the formula
\begin{equation}\label{eq:beta-0}
B(n)=-\S(n)-\sum_{s=1}^{n-1}B(s)\S(n-s)\, , 
\end{equation}
then 
\begin{equation}\label{eq:beta-1}
\S^+(n)=\S(n)+\sum_{s=1}^{R-1}B(s)\S(n-s)+O(\S(n-R))
\end{equation}
for each integer $R>1$. In particular, axiom $\cG_3^+$ is also satisfied.
\end{lem}

\medskip\noindent\emph{Proof:} According to \cite{K}, the sequences $\{\S(n)\}$ and $\{\S^+(n)\}$ are related by the additive Euler product formula
\[
\sum_{n=0}^\infty \S(n) x^n =\prod_{m=1}^\infty \left( 1-x^m \right)^{-\S^+(m)}
\]
understood as an equality of formal power series. Since by assumption $\S(n)\ge \S^+(n)>0$ for $n$ sufficiently large, the result follows by Theorem 5 and Theorem 7 in \cite{W3}. \qed

\begin{ex}\label{ex:graph-counting}
Let $(\G,+,\Box,\partial)$ be the arithmetical semiring of graphs introduced in Example \ref{ex:graph-semiring}. As shown in \cite{W4}, the sequence $\{\G(n)\}$ has asymptotic expansion
\begin{equation}\label{eq:numberofgraphs}
\G(n)= \frac{2^{\binom{n}{2}}}{n!}+\sum_{s=1}^{R-1}\varphi_s(n)\frac{2^{\binom{n-s}{2}}}{(n-s)!}+O\left(\frac{2^{\binom{n-R}{2}}}{(n-2R)!}\right)
\end{equation}
where $R$ is any integer greater than $1$ and, for each $s\in \N$, $\varphi_s(n)$ is polynomial of degree $s$ in $n$. An explicit group-theoretic formula for $\varphi_s$ is given in \cite{W4}. For small values of $s$ it yields
\begin{eqnarray*}
\varphi_1(n)&=&n-1\,;\\
\varphi_2(n)&=&\frac{1}{3}(3n^2-13n+14)\,;\\
\varphi_3(n)&=&\frac{1}{3}(4n^3-46n^2+177n-225)\,;\\
\varphi_4(n)&=&\frac{2}{45}(60n^4-1305n^3+10580n^2-37737n+49828)\,.
\end{eqnarray*}
Therefore \eqref{eq:log} holds with $\alpha=\log(\sqrt{2})$, $\beta=-1$, $\gamma=-\alpha-\beta$, $a=1$ and, say, $b=\frac{1}{2}$. By Lemma \ref{lem:log}, the arithmetical semiring $(\G,+,\Box,\partial)$ satisfies axioms $\cG_1^+$, $\cG_2^+$ and $\cG_3^+$. Substituting \eqref{eq:numberofgraphs} into \eqref{eq:beta-1} one obtains the asymptotic expansion 
\begin{equation}\label{eq:numberofconnectedgraphs}
\G^+(n)=\frac{2^{\binom{n}{2}}}{n!}+\sum_{s=1}^{R-1}\omega_s(n)\frac{2^{\binom{n-s}{2}}}{(n-s)!} + O\left(\frac{2^{\binom{n-R}{2}}}{(n-2R)!}\right)
\end{equation}
found in \cite{W3}, where $R$ is any integer greater than $1$ and for all $s\in\N$
\[
\omega_s(n)=\varphi_s(n)+B(s)+\sum_{r=1}^{s-1}B(r)\varphi_s(n-r) 
\]
where $B$ is defined by specializing \eqref{eq:beta-0} to $\S=\G$. For instance, using $\G(1)=1$, $\G(2)=2$, $\G(3)=4$ and $\G(4)=11$ one obtains $B(1)=B(2)=B(3)=-1$, $B(4)=-4$ and thus \cite{W3} 
\begin{eqnarray*}
\omega_1(n)&=&n-2\,;\\
\omega_2(n)&=&\frac{1}{3}(3n^2-16n+17)\,;\\
\omega_3(n)&=&\frac{1}{3}(4n^3-49n^2+193n-249)\,;\\
\omega_4(n)&=&\frac{1}{45}(120n^4-2670n^3+21985n^2-79359n+105656)\,.
\end{eqnarray*}
\end{ex}

\begin{ex}\label{ex:evengraphs}
As observed in \cite{liskovets}, $2\G_{\rm ev}(n)-\G(n)$ equals the number of self-complementary graphs with $n$ vertices. Therefore, using the asymptotic estimates of \cite{palmer}, we obtain
\[
\log(2\G_{\rm ev}(n)-\G(n))=o\left(\log\left(\G\left(3\left\lfloor\frac{n}{2}\right\rfloor\right)\right)\right)
\]
and thus
\[
2\G_{\rm ev}(n)-\G(n)=O\left(\G\left(3\left\lfloor\frac{n}{2}\right\rfloor\right)\right)=O\left(\frac{2^{\binom{n-R}{2}}}{(n-2R)!}\right)
\]
for all integers $R>1$. Therefore
\[
\G_{\rm ev}(n)= \frac{2^{\binom{n}{2}}}{2n!}+\sum_{s=1}^{R-1}\varphi_s(n)\frac{2^{\binom{n-s}{2}}}{2(n-s)!}+O\left(\frac{2^{\binom{n-R}{2}}}{(n-2R)!}\right)
\]
where the polynomials $\varphi_s$ are as in \eqref{eq:numberofgraphs} and $R$ is any integer greater than $1$. In particular, the arithmetical semiring $(\G_{\rm ev},+,\Box,\partial)$ satisfies the assumptions of Lemma \ref{lem:log} and thus axioms $\cG_1^+$, $\cG_2^+$ and $\cG_3^+$. Moreover, for any integer $R>1$,
\begin{equation}\label{eq:numberofevenconnectedgraphs}
\G_{\rm ev}^+(n)=\frac{2^{\binom{n}{2}}}{2n!}+\sum_{s=1}^{R-1}\omega_s^{\rm ev}(n)\frac{2^{\binom{n-s}{2}}}{2(n-s)!} + O\left(\frac{2^{\binom{n-R}{2}}}{(n-2R)!}\right)\,,
\end{equation}
where
\[
\omega_s^{\rm ev}(n)=\varphi_s(n)+B_{\rm ev}(s)+\sum_{r=1}^{s-1}B_{\rm ev}(r)\varphi_s(n-r) 
\]
for each $s\in \N$ and $B_{\rm ev}$ is the sequence obtained by specializing \eqref{eq:beta-0} to $\S=\G_{\rm ev}$. For instance, using $\G_{\rm ev}(1)=1$, $\G_{\rm ev}(2)=1$, $\G_{\rm ev}(3)=2$, $\G_{\rm ev}(4)=6$ we obtain $B(1)=1$, $B(2)=0$, $B(3)=-1$ and $B(4)=-3$ and thus
\begin{eqnarray*}
\omega_1^{\rm ev}(n)&=&n-2\,;\\
\omega_2^{\rm ev}(n)&=&\frac{1}{3}(3n^2-16n+20)\,;\\
\omega_3^{\rm ev}(n)&=&\frac{1}{3}(4n^3-49n^2+196n-258)\,;\\
\omega_4^{\rm ev}(n)&=&\frac{1}{45}(120n^4-2670n^3+22030n^2-79734n+106481)\,.
\end{eqnarray*}
\end{ex}

\begin{lem}\label{lem:1}
Let $(\S,+,\Box,\partial)$ be an arithmetical semiring. Then for each $n\in \mathbb N$,
\begin{equation}\label{eq:1}
\sum_{r=1}^{\lceil n/2\rceil-1}\S^+(r)\S^+(n-r)\le \S(n)-\S^+(n)-\binom{\S^+(n/2)}{2}-\S^+\left(\frac{n}{2}\right)\le \sum_{r=1}^{\lceil n/2\rceil-1}\S^+(r)\S(n-r)
\end{equation}
and
\begin{equation}\label{eq:2}
\sum_{r=2}^{\lceil \sqrt{n} \rceil-1} \S^\Box(r)\S^\Box(n/r) \le \S^+(n)-\S^\Box(n)-\binom{\S^{\Box}(\sqrt{n})}{2}-\S^\Box(\sqrt{2})\le\sum_{r=2}^{\lceil \sqrt{n} \rceil-1}  \S^\Box(r)\S^+(n/r)\,.
\end{equation}
Here the notation of Remark \ref{rem:3} is tacitly employed so that, in particular, the summation indices $r$ in \eqref{eq:2} can be equivalently thought of as being further restricted to divisors of $n$.   
\end{lem}

\medskip\noindent\emph{Proof:} Fix a total order $\prec$ on $\S$ in such a way that $S'\prec S''$ whenever $\partial(S')<\partial(S'')$ and let
\[
F_+:\S\setminus \{\mathbf e_+\}\to \S^+
\]
be the function that to each $S\neq \mathbf e_+$ assigns the $\prec$-smallest additive prime in the unique additive factorization of $S$. By unique factorization, we obtain a second function
\[
F'_+:\S\setminus \{\mathbf e_+\}\to \S
\]
defined by declaring $F'_+(S)$ to be the unique element such that $S=F_+(S)+F'_+(S)$. For each $n\ge 1$ and for each $r=1,\ldots,\lfloor n/2\rfloor$, let 
\[
X_{n,r}^+=\{S\in \uS(n)\setminus \uS^+(n)\,|\,\partial(F_+(S))=r\}\,.
\]
Since $X_{n,r}^+\cap X_{n,r'}^+=\emptyset$ whenever $r\neq r'$, then
\begin{equation}\label{eq:lem24-1}
\S(n)-\S^+(n)=\left|\bigcup_{r=1}^{\lfloor n/2\rfloor} X_{n,r}^+  \right| = \sum_{r=1}^{\lfloor n/2\rfloor }\left|X_{n,r}^+\right|\,.
\end{equation}
By construction, $S\in X_{n,n/2}^+$ if and only if both $F_+(S)$ and $F'_+(S)$ are additive primes such that $\partial(F_+(S))=\partial(F'_+(S))=\frac{n}{2}$. Therefore,
\begin{equation}\label{eq:lem24-2}
\left| X_{n,n/2}^+\right| = \binom{\S^+(n/2)}{2}+\S^+\left(\frac{n}{2}\right)\,,
\end{equation}
where the second term accounts for additive squares (i.e.\ doubles). On the other hand, for each $r=1,\ldots,\lceil n/2\rceil -1$ there are injections 
\[
\uS^+(r)\times \uS^+(n-r)\to X_{n,r}^+ \to \uS^+(r)\times \uS(n-r)
\]
where the first map is addition and the second map is the assignment $S\mapsto (F_+(S),F'_+(S))$. Hence
\begin{equation}\label{eq:lem24-3}
\sum_{r=1}^{\lceil n/2\rceil-1} \S^+(r)\S^+(n-r)\le \sum_{n=1}^{\lceil n/2\rceil-1} \left| X_{n,r}^+ \right| \le \sum_{r=1}^{\lceil n/2\rceil-1} \S^+(r)\S(n-r)
\end{equation}
from which, using \eqref{eq:lem24-1} and \eqref{eq:lem24-2}, \eqref{eq:1} easily follows. Similarly, let $F_\Box:\S^+\setminus \{\mathbf e_\Box\}\to \S^\Box$ and $F_\Box':\S^+\setminus \{\mathbf e_\Box\}\to \S^+$ such that for each additive prime $S$, $F_\Box(S)$ is the $\prec$-smallest multiplicative prime divisor of $S$ and $S=F_\Box(S)\Box F_\Box'(S)$. If 
\begin{equation}\label{eq:lem24-4}
X_{n,r}^\Box=\{S\in \uS^+(n)\setminus \uS^\Box(n)\,|\,\partial(F_\Box(S))=r\}\,
\end{equation}
then 
\[
\S(n)-\S^\Box(n)= \sum_{r=2}^{\lfloor \sqrt{n}\rfloor} \left|X_{n,r}^\Box\right| = \binom{\S^\Box(\sqrt{n})}{2}+\S^\Box(\sqrt{n})+\sum_{r=2}^{\lceil\sqrt{n}\rceil-1} \left|X_{n,r}^\Box\right|
\]
and thus \eqref{eq:2} follows from the existence of injections
\begin{equation}\label{eq:lem24-5}
\uS^\Box(r)\times \uS^\Box\left(\frac{n}{r}\right)\to X_{n,r}^\Box \to \uS^\Box(r)\times \uS^+\left(\frac{n}{r}\right)
\end{equation}
where the first map is multiplication and the second is the assignment $S\mapsto (F_\Box(S),F'_\Box(S))$.
\qed

\begin{lem}\label{lem:2}
Let $(\S,+,\Box,\partial)$ be a monotonic arithmetical semiring and let $D>2$ be an integer. Then 
\begin{equation}\label{eq:lem2}
0\le \S^+(n)-\S^\Box(n) - \sum_{r=2}^{D-1} \S^\Box(r)\S^+\left(\frac{n}{r}\right) \le \S\left(\left\lfloor\frac{n}{D}\right\rfloor +D \right)-\S^+\left(\left\lfloor\frac{n}{D}\right\rfloor +D \right)
\end{equation}
for all sufficiently large $n$.
\end{lem}

\medskip\noindent\emph{Proof:}
If $r$ is a divisor of $n$ such that  $D\le r\le \lfloor n/D \rfloor$, then $n\ge rD$. Hence $n\left(\frac{1}{D}-\frac{1}{r} \right)\ge (r-D)$, which implies
\begin{equation}\label{eq:3}
\frac{n}{r}\le \left\lfloor \frac{n}{D}\right\rfloor+D-r\,.
\end{equation} 
It follows from \eqref{eq:2} that
\begin{equation}\label{eq:4}
\S^+(n)-\S^\Box(n)\le \binom{\S^\Box(\sqrt{n})}{2}+ \S^\Box(\sqrt{n})+\sum_{r=2}^{D-1}\S^\Box(r)\S^+\left(\frac{n}{r}\right) + \sum_{r=D}^{\lceil \sqrt{n} \rceil -1} \S^\Box(r)\S^+\left(\frac{n}{r}\right)
\end{equation}
Since $\sqrt{n}\le \left\lfloor \frac{n}{D}\right\rfloor+D-\sqrt{n}$ for all $n$ sufficiently large, using monotonicity yields
\begin{equation}\label{eq:4.1}
\binom{\S^\Box(\sqrt{n})}{2}+\S^\Box(\sqrt{n})\le \S^\Box(\sqrt{n}) \S^\Box(\sqrt{n})\le \S^+(\sqrt{n})\S^+\left(\left\lfloor \frac{n}{D}\right\rfloor+D-\sqrt{n}\right)\,.
\end{equation}
Using monotonicity in a similar way results in the upper bound
\begin{equation}\label{eq:4.2}
\sum_{r=D}^{\lceil \sqrt{n} \rceil -1} \S^\Box(r)\S^+\left(\frac{n}{r}\right)  \le 
\sum_{r=D}^{\lceil\sqrt{n}\rceil-1} \S^+(r)\S^+\left(\left\lfloor \frac{n}{D}\right\rfloor+D-r\right)
\end{equation}
which is valid for all $n$ sufficiently large. Combining \eqref{eq:4.1} with \eqref{eq:4.2} and using \eqref{eq:1} we obtain the estimate
\begin{eqnarray*}
\binom{\S^\Box(\sqrt{n})}{2}+\S^\Box(\sqrt{n})+\sum_{r=D}^{\lceil \sqrt{n} \rceil -1} \S^\Box(r)\S^+\left(\frac{n}{r}\right) &\le& \sum_{r=D}^{\lfloor \sqrt{n} \rfloor} \S^+(r)\S^+\left(\left\lfloor \frac{n}{D}\right\rfloor+D-r\right)\\
&\le & \S\left(\left\lfloor\frac{n}{D}\right\rfloor +D \right)-\S^+\left(\left\lfloor\frac{n}{D}\right\rfloor +D \right)
\end{eqnarray*}
which, substituted into \eqref{eq:4}, yields the second inequality in \eqref{eq:lem2}. On the other hand, monotonicity implies
\[
\sum_{r=2}^{D-1}\S^\Box(r)\S^+\left(\frac{n}{r}\right) \le \sum_{r=2}^{D-1}\S^+(r)\S^+(n-r)
\]
and thus, using \eqref{eq:1}, the first inequality in \eqref{eq:lem2}. 
\qed

\begin{thm}\label{thm:1}
Let $(\S,+,\Box,\partial)$ be a monotonic arithmetical semiring and let $p$ be the smallest integer such that such that $\S^\Box(p)\neq 0$. Then
\begin{equation}\label{eq:thm}
\S^+(n)-\S^\Box(n)= \S^\Box(p)\S^+\left(\frac{n}{p}\right)+O\left(\S\left(\left\lfloor \frac{n}{p+1}\right\rfloor+p+1\right)\right)\,.
\end{equation}
Moreover, if axiom $\cG_1^+$ holds then axiom $\cG_1^\Box$ also holds. 
\end{thm}

\medskip\noindent\emph{Proof:}
Setting $D=p+1$ in Lemma \ref{lem:2} we obtain
\begin{equation}\label{eq:thm26.1}
\S^+(n)-\S^\Box(n)-\S^\Box(p)\S^+\left(\frac{n}{p}\right) \le \S\left(\left\lfloor\frac{n}{p+1}\right\rfloor +p+1 \right)-\S^+\left(\left\lfloor\frac{n}{p+1}\right\rfloor +p+1 \right)
\end{equation}
from which \eqref{eq:thm} easily follows. 
On the other hand, $\cG_1^+$ implies 
\begin{equation}\label{eq:thm26.2}
\S\left(\left\lfloor\frac{n}{p+1}\right\rfloor +p+1 \right)-\S^+\left(\left\lfloor\frac{n}{p+1}\right\rfloor +p+1 \right)=o\left(\S^+\left(\left\lfloor\frac{n}{p+1}\right\rfloor +p+1 \right) \right)\,.
\end{equation}
It follows from monotonicity and \eqref{eq:thm26.2} that the right hand side of \eqref{eq:thm26.1} is $o(\S^+(n))$. Moreover, $\cG_1^+$ also implies
\begin{equation}\label{eq:thm26.3}
\S^\Box(p)\S^+\left(\frac{n}{p}\right)\le \S(n)-\S^+(n)=o(\S^+(n))\,.
\end{equation}
where the inequality is a consequence of \eqref{eq:lem2}.
Substitution into \eqref{eq:thm26.1} yields $\S^+(n)-\S^\Box(n)=o(\S^+(n))$, which is equivalent to the validity of axiom $\cG_1^\Box$.
\qed

\begin{cor}\label{cor:1}
If $(\S,+,\Box,\partial)$ is an arithmetical semiring that satisfies axioms $\cG_1^+$ and $\cG_2^+$ then it also satisfies axioms $\cG_1^\Box$, $\cG_2^\Box$ and $\cG_3^\Box$.
\end{cor}

\medskip\noindent\emph{Proof:} By Remark \ref{rem:monotonic}, $(\S,+,\Box,\partial)$ is monotonic and thus, by Theorem \ref{thm:1}, axiom $\cG_1^\Box$ holds. On the other hand, 
\begin{equation}\label{eq:13}
\S^\Box\left(\left\lfloor \frac{n}{2}\right\rfloor\right)\le \S^+\left(\left\lfloor \frac{n}{2}\right\rfloor\right) = o(\S^+(n))\,
\end{equation}
where the equality follows from iterated use of axiom $\cG_2^+$. Combining \eqref{eq:13} with axiom $\cG_1^\Box$ yields axiom $\cG_2^\Box$. Moreover, if $p$ is the smallest integer for which $\S^\Box(p)\neq 0$, then
\begin{equation}\label{eq:13.1}
\S\left(\left\lfloor \frac{n}{p+1}\right\rfloor+p+1\right) = O\left(\S^+\left(\left\lfloor \frac{n}{p+1}\right\rfloor+p+1\right)\right) 
\end{equation}
by axiom $\cG_1^+$. Iterated use of axiom $\cG_2^+$ shows that the right hand side of \eqref{eq:13.1} is $o\left(\S^+\left(\left\lfloor \frac{n}{2}\right\rfloor\right)\right)$ and thus $o\left(\S^\Box\left(\left\lfloor \frac{n}{2}\right\rfloor\right)\right)$  by axiom $\cG_1^\Box$. A further application of axiom $\cG_1^\Box$ yields
\begin{equation}\label{eq:13.2}
\S^\Box(p)\S^+\left(\frac{n}{p} \right)\le \S^\Box(p)\S^+ \left(\left\lfloor \frac{n}{2}\right\rfloor \right) = O\left(\S^\Box \left(\left\lfloor \frac{n}{2}\right\rfloor \right)\right)
\end{equation}
which, upon substitution into the statement of Theorem \ref{thm:1}, completes the proof that $\cG_3^\Box$ holds. 
\qed

\begin{cor}\label{cor:2} Let $(\S,+,\Box,\partial)$ be a monotonic arithmetical semiring and let $p$ be the smallest integer such that $\S^\Box(p)\neq 0$. Then
\begin{equation}\label{eq:13.9}
\S^+(n)-\S^\Box(n)\le \S^\Box(p)\S\left(\left\lfloor \frac{n}{p} \right\rfloor \right) + \S\left(\left\lfloor \frac{n}{p+1} \right\rfloor+p+1 \right)
\end{equation}
for all $n$ sufficiently large. Moreover, if axiom $\cG_1^+$ holds, then
\begin{equation}\label{eq:14}
\S^+(pn)-\S^\Box(pn)=\S^\Box(p)\S^+(n)+o\left(\S^+\left(\left\lfloor\frac{pn}{p+1}\right\rfloor+p+1\right)\right)\,.
\end{equation}
\end{cor}

\medskip\noindent\emph{Proof:} Substitution of
\[
\S^+\left(\frac{n}{p}\right) \le \S\left(\left\lfloor\frac{n}{p}\right\rfloor\right)\
\]
into \eqref{eq:thm26.1} yields \eqref{eq:13.9}. The second statement follows upon substitution of  \eqref{eq:thm26.2} into \eqref{eq:thm26.1}.
\qed

\begin{rem}\label{rem:log} Let $(\S,+,\Box,\partial)$ be an arithmetical semiring such that $\eqref{eq:log}$ holds, let $p$ be the smallest integer such that $\S^\Box(p)\neq 0$ and let $B$ be defined by \eqref{eq:beta-0}.Then combining \eqref{eq:beta-1} and \eqref{eq:14} we obtain the asymptotic expansion
\begin{equation}\label{eq:difference-log}
\S^+(pn)-\S^\Box(pn)=\S^\Box(p)\left(\S(n) +\sum_{s=1}^{R-1} B(s) \S(n-s)\right)+O(\S(n-R))\,,
\end{equation}
for every integer $R\ge 2$. 
\end{rem}

\begin{ex}
Consider our main example, the arithmetical semiring of graphs $(\G,+,\Box,\partial)$. By Example \ref{ex:graph-counting}, we may apply Lemma \ref{lem:log} to conclude the validity of axioms $\cG_1^+$, $\cG_2^+$ and $\cG_3^+$. By Corollary \ref{cor:1}, we conclude that axioms $\cG_1^\Box$, $\cG_2^\Box$ and $\cG_3^\Box$ also hold. In particular
\[
\G(n)\sim \G^+(n)\sim \G^\Box(n)
\]
i.e.\ almost all graphs are both connected and cartesian prime. Since $\underline \G^\Box(2)=\{K_2\}$, then Corollary \ref{cor:2} and Remark \ref{rem:log} imply the more precise estimates \eqref{eq:0} and \eqref{eq:00}, respectively. 
\end{ex}

\begin{ex} Consider now the arithmetical semiring $(\G_{\rm ev},+,\Box,\partial)$ of graphs with an even number of edges. By Example \ref{ex:evengraphs}, we may apply Lemma \ref{lem:log} to conclude the validity of axioms $\cG_1^+$, $\cG_2^+$ and $\cG_3^+$. By Corollary \ref{cor:1}, we conclude that axioms $\cG_1^\Box$, $\cG_2^\Box$ and $\cG_3^\Box$ also hold. In particular
\[
\G_{\rm ev}(n)\sim \G^+_{\rm ev}(n)\sim \G_{\rm ev}^\Box(n)
\]
i.e.\ almost all graphs with an even number of edges are both connected and prime with respect to the product induced by the cartesian product. Since the only connected graph with two vertices has exactly one edge, then $\G_{\rm ev}^\Box(2)=0$. On the other hand, $\G_{\rm ev}^\Box(3)=1$ and thus 
\[
\G_{\rm ev}(n)-\G_{\rm ev}^\Box \le \G_{\rm ev}\left(\left\lfloor \frac{n}{3}\right\rfloor \right) + \G_{\rm ev}\left(\left\lfloor \frac{n}{4}\right\rfloor+4 \right) 
\]
by Corollary \ref{cor:2}. Moreover, Remark \ref{rem:log} implies that for each integer $R>1$
\[
\G_{\rm ev}^+(3n)-\G_{\rm ev}^\Box(3n)=\frac{2^{\binom{n}{2}}}{2n!}+\sum_{s=1}^{R-1}\omega_s^{\rm ev}(n)\frac{2^{\binom{n-s}{2}}}{2(n-s)!} + O\left(\frac{2^{\binom{n-R}{2}}}{(n-2R)!}\right)
\]
where $\{\omega_s^{\rm ev}\}$ is the sequence of polynomials defined in Example \ref{ex:evengraphs}.

\end{ex}

\section{submultiplicative functions}

In this section we derive additional information about the distribution of multiplicative primes by looking at real valued functions of additive primes that are compatible with multiplication. Our main result states that, under natural growths for the number of additive primes, these functions are asymptotically dominated by their restriction to multiplicative primes. The examples that we choose to illustrate our results are direct multiplicative analogues of the function discussed in \cite{K} in the context of additive arithmetical semigroups.    

\begin{dfn} Let $(\S,+,\Box,\partial)$ be an arithmetical semiring. A function $f:\S^+\to \R_{\ge 0}$ is {\it submultiplicative} if  $f(S_1\Box S_2)\le f(S_1)f(S_2)$ for all $S_1,S_2\in \S^+$. Given $f$ submultiplicative, and $\bullet\in\{+,\Box\}$ we define sequences 
$f^\bullet$, $f^\bullet_{\rm max}$ such that
\[
f^\bullet(n)=\sum_{S\in \uS^\bullet(n)} f(S)\quad\textrm{ and }\quad f^\bullet_{\rm max}(n)=\max_{s\le n}\max_{S\in \uS^\bullet(n)}f(S)
\]
for all $n\in \N$. We say that $f$ has {\it asymptotic-mean value} $\mu$ on $\S^\bullet$ if
\[
\lim_{n\to \infty} \frac{f^\bullet(n)}{\S^\bullet(n)} = \mu\,.
\]
If this is the case, we say that $f$ has asymptotic variance $\nu$ on $\S^\bullet$ if
\[
\lim_{n\to \infty} \frac{1}{\S^\bullet(n)} \sum_{S\in \uS^\bullet(n)}(f(S)-\mu)^2=\nu\,.
\]
\end{dfn}

\begin{lem}\label{lem:1.9}
Let $f$ be a submultiplicative function on an arithmetical semiring $(\S,+,\Box,\partial)$. Then
\[
f^+(n)-f^\Box(n)\le \left(f_{\rm max}^\Box(\sqrt{n})\right)^2\left( \binom{\S^\Box(\sqrt{n})}{2}+\S^\Box(\sqrt{n})\right)+\sum_{r=2}^{\lceil \sqrt{n} \rceil -1} f_{\rm max}^\Box(r) f_{\rm max}^+\left(\frac{n}{r}\right)\S^\Box(r)\S^+\left(\frac{n}{r}\right)\,.
\]
\end{lem}

\medskip\noindent\emph{Proof:} Using the partition
\[
\uS^+(n)\setminus \uS^\Box(n) = \bigcup_{r=2}^{\lfloor \sqrt{n}\rfloor} X_{n,r}^\Box
\]
and the functions $F_\Box$, $F_\Box'$ introduced in the proof of Lemma \ref{lem:1}, we obtain
\begin{equation}\label{eq:lem191}
f^+(n)-f^\Box(n)=\sum_{r=2}^{\lfloor \sqrt{n}\rfloor} \sum_{S\in X_{n,r}^\Box} f(S)\le \sum_{r=2}^{\lfloor \sqrt{n}\rfloor} \sum_{S\in X_{n,r}^\Box} f(F_\Box(S))f(F_\Box'(S)) \,.
\end{equation}
If $r<\sqrt{n}$
\begin{equation}\label{eq:lem192}
\sum_{S\in X_{n,r}^\Box} f(F_\Box(S))f(F_\Box'(S)) \le f_{\rm max}^\Box(r) f_{\rm max}^+\left(\frac{n}{r}\right)\left|X_{n,r}^\Box \right| \le f_{\rm max}^\Box(r) f_{\rm max}^+\left(\frac{n}{r}\right)\S^\Box(r)\S^+\left(\frac{n}{r}\right)\,,
\end{equation}
where the second inequality follows since the second map in \eqref{eq:lem24-5} is an injection.
On the other hand, since $X_{n,\sqrt n}^\Box$ is in bijection with the set of unordered pairs of multiplicative primes in $\uS^\Box(\sqrt n)$ then
\begin{equation}\label{eq:lem193}
\sum_{S\in X_{n,\sqrt{n}}^\Box} f(F_\Box(S))f(F_\Box'(S)) \le \left(f_{\rm max}^\Box(\sqrt{n})\right)^2\left( \binom{\S^\Box(\sqrt{n})}{2}+\S^\Box(\sqrt{n})\right)\,.
\end{equation}
Substitution of \eqref{eq:lem193} and \eqref{eq:lem192} into \eqref{eq:lem191} then proves the lemma.
\qed

\begin{thm}\label{thm:2}
Let $f$ be a submultiplicative function on a monotonic arithmetical semiring $(\S,+,\Box,\partial)$ such that $f^+_{\rm max}(n)=O(\S^+(n))$. Then
\begin{equation}\label{eq:function-0}
f^+(n)-f^\Box(n)=O(\S(n)-\S^+(n))\,. 
\end{equation}
\end{thm}

\noindent\emph{Proof:} By assumption there exists positive integers $N, C$ such that
\[
f_{\rm max}^+(n)\le C \S^+(n)\le C \S^+(n+1)
\]
for all $n\ge N$. Then for $r\ge N$
\begin{equation*}\label{eq:thm2-1}
f_{\rm max}^\Box(r)f_{\rm max}^+\left(\frac{n}{r}\right) \S^\Box(r)\S^+\left(\frac{n}{r} \right) \le C^2 P_{n,r}
\end{equation*}
where
\[
P_{n,r}=\S^+(r)\S^+(r+1)\S^+\left( \frac{n}{r}+1\right)\S^+ \left(n-\frac{n}{r}-2r-2\right).
\]
If $2<r<\sqrt{n}$, then
\[
r<r+1<\frac{n}{r}+1<n-\frac{n}{r}-2r-2
\]
and thus $P_{n,r}$ is the number of elements $S\in X^+_{n,r}$ such that $F_+'(S)=S_1+ S_2+ S_3$
for some $S_1\in \uS^+(r+1)$, $S_2\in \uS^+\left(\frac{n}{r}+1 \right)$, and $S_3\in \uS^+\left(n-\frac{n}{r}-2r-2 \right)$. Since $X^+_{n,r}\cap X^+_{n,r'}=\emptyset$ whenever $r\neq r'$, then
\[
\sum_{r=N}^{\sqrt{n}\rceil-1} P_{n,r}\le \sum_{r=N}^{\lceil \sqrt{n}\rceil-1} \left|X^+_{n,r}\right| = \left|\bigcup_{r=N}^{\lceil \sqrt{n}\rceil-1} X^+_{n,r}\right| \le \S(n)-\S^+(n)\,. 
\]
and thus
\begin{equation}\label{eq:thm2-2}
\sum_{r=N}^{\lceil \sqrt{n} \rceil -1}  f_{\rm max}^\Box(r)f_{\rm max}^+\left(\frac{n}{r}\right) \S^\Box(r)\S^+\left(\frac{n}{r} \right) = O(\S(n)-\S^+(n))\,. 
\end{equation}
Similarly,
\begin{align}
\left(f_{\rm max}^\Box(\sqrt{n})\right)^2\frac{(\S^\Box(\sqrt{n}))^2+\S^\Box(\sqrt{n})}{2} & \le  C^2 (\S^+(\sqrt{n}))^4\nonumber\\
&\le 
C^2 \S^+(\sqrt{n})\S^+(\sqrt{n}+1)\S^+(\sqrt{n}+2)\S^+(n-3\sqrt{n}-3)\nonumber\\  & \le  C^2(\S(n)-\S^+(n))\,.\label{eq:thm2-3}
\end{align}
Moreover
\begin{equation}\label{eq:thm2-4}
\sum_{r=2}^{N-1} f_{\rm max}^\Box(r)f_{\rm max}^+\left(\frac{n}{r}\right) \S^\Box(r)\S^+\left(\frac{n}{r} \right) = O\left(\S^+\left(\frac{n}{p}\right)\S^+\left(\frac{n}{p}\right)\right) 
\end{equation}
where $p$ is the smallest integer such that $\S^\Box(p)\neq 0$. If $p>2$, then
\[
\S^+\left(\frac{n}{p}\right)\S^+\left(\frac{n}{p}\right)\le \S^+\left(\frac{n}{p}\right)\S^+\left(n-\frac{n}{p}\right) \le \S(n)-\S^+(n)\,.
\]
On the other hand 
\begin{equation}\label{eq:41}
\S^+\left(\frac{n}{2}\right)\S^+\left(\frac{n}{2}\right)\le \binom{\S^+(n/2)}{2} + \S^+\left(\frac{n}{2}\right)\le \S(n)-\S^+(n)
\end{equation}
and thus, using \eqref{eq:thm2-4}, we obtain 
\begin{equation}\label{eq:thm2-5}
\sum_{r=2}^{N-1} f_{\rm max}^\Box(r)f_{\rm max}^+\left(\frac{n}{r}\right) \S^\Box(r)\S^+\left(\frac{n}{r} \right) = O(\S(n)-\S^+(n))
\end{equation}
independently of $p$. Finally, substitution of \eqref{eq:thm2-2}, \eqref{eq:thm2-3}, and \eqref{eq:thm2-5} into the statement of Lemma \ref{lem:1.9} yields \eqref{eq:function-0}. \qed

\begin{cor}\label{cor:3}
Let $(\S,+,\Box,\partial)$ be a monotonic arithmetical semiring that satisfies axiom $\cG_1^+$ and let $f$ be a submultiplicative function on $(\S,+,\Box,\partial)$ that has asymptotic mean-value $\mu$ on $\S^\Box$. Then
\begin{enumerate}[1)]
\item $f_{\rm max}^+(n)=O(\S^+(n))$ implies that $f$ has asymptotic mean-value $\mu$ on $\S^+$; 
\item $f_{\rm max}^+(n)=O\left(\S^+\left(\left\lfloor\frac{n}{2}\right\rfloor\right)\right)$ implies that $f$ has asymptotic variance $0$ on $\S^+$ if it has asymptotic variance $0$ on $\S^\Box$. 
\end{enumerate}
\end{cor}
 
\noindent\emph{Proof:} Axiom $\cG_1^+$ applied to \eqref{eq:function-0} implies $f^+(n)-f^\Box(n)=o(\S^+(n))$. Therefore, using axiom $\cG_1^+$ once more, we obtain
\[
\mu=\lim_{n\to \infty} \frac{f^\Box(n)}{\S^\Box(n)}=\lim_{n\to \infty} \frac{f^\Box(n)}{\S^+(n)}=\lim_{n\to \infty} \frac{f^+(n)}{\S^+(n)}\,.
\]
If
$f_{\rm max}^+(n)=O\left(\S^+\left(\left\lfloor\frac{n}{2}\right\rfloor\right)\right)$, then 
\[
(f^2)_{\rm max}^+(n)\le (f_{\rm max}^+(n))^2 = O\left(\S^+\left(\left\lfloor\frac{n}{2}\right\rfloor\right)\S^+\left(\left\lfloor\frac{n}{2}\right\rfloor\right)\right)
\]
from which, using \eqref{eq:41} and axiom $\cG_1^+$, we obtain $(f^2)_{\rm max}^+(n)=O(\S^+(n))$. By assumption
\begin{equation}\label{eq:variance}
0=\lim_{n\to \infty}\frac{1}{\S^\bullet(n)}\sum_{S\in \uS^\bullet(n)}(f(S)-\mu)^2=\lim_{n\to \infty}\frac{(f^2)^\bullet(n)-2\mu f^\bullet(n)+\mu^2\S^\bullet(n)}{\S^\bullet(n)}
\end{equation}
holds for $\bullet=\Box$. Hence, $f^2$ has asymptotic mean-value $\mu^2$ on $\S^\Box$ and thus, by Theorem \ref{thm:2}, on $\S^+$. This implies that \eqref{eq:variance} holds for $\bullet=+$, which proves 2). \qed

\begin{ex}\label{ex:divisor}
Let $(\S,+,\Box,\partial)$ be a monotonic arithmetical semiring for which $\cG_1^+$ holds. The {\it divisor function} $d:\S^+\to \R_{\ge 0}$ is the function that to each $S\in \S^+$ assigns the total number $d(S)$ of factorizations of the form $S=S'\Box S''$ for some $S',S''\in \S^+$. Since $d(S)=2$ for all $S\in \S^\Box$ and the total number of non-trivial divisors of $S\in\uS^+(n)$ is at most the total number of additive primes whose degree divides $n$, then 
\begin{eqnarray*}
d^+_{\rm max}(n)&\le&2+\sum_{r=2}^{\lfloor n/2\rfloor}\S^+(r) \\
&\le& 2 + \S^+(2)+\S^+\left(\left\lfloor\frac{n}{2}\right\rfloor\right) + \sum_{r=3}^{\lfloor n/3\rfloor}\S^+(r)\S^+\left(\left\lfloor\frac{n}{3}\right\rfloor+3-r \right)\\
&=& 2+ \S^+(2) + O\left(\S^+\left(\left\lfloor\frac{n}{2}\right\rfloor\right) \right)
\end{eqnarray*}
where \eqref{eq:1}, axiom $\cG_1^+$, and monotonicity were used to derive the last equality. If we additionally assume that the sequence $\{\S^+(n)\}$ is unbounded (which is the case if the arithmetical semiring is strictly monotone), then $d^+_{\rm max}(n)= O\left(\S^+\left(\left\lfloor\frac{n}{2}\right\rfloor\right) \right)$ and thus $d$ satisfy the assumptions of Corollary \ref{cor:3}. We conclude that $d$ has asymptotic mean-value $2$ and asymptotic variance $0$ on $\S^+$. 
\end{ex}

\begin{ex}\label{ex:unitary-divisor}
Let $(\S,+,\Box,\partial)$ be a strictly monotonic arithmetical semiring for which axiom $\cG_1^+$ holds. Let $d_*:\S^+\to \R_{\ge 0}$ be the {\it unitary-divisor function} that to each $S\in \S^+$ assigns the number $d_*(S)$ of factorizations $S=S'\Box S''$ for some $S',S''\in \S^+$ that are {\it $\Box$-coprime} i.e.\ such that $S'$ and $S''$ have no common $\Box$-prime factors. Since $d_*(S)\le d(S)$ for each $S\in \S^+$ and thus $(d_*)^+_{\rm max}\le d^+_{\rm max}$, then the argument of Example \ref{ex:divisor} applies to $d_*$ as well. In particular $d_*$ has asymptotic mean value $2$ and asymptotic variance $0$. Similarly, let $\beta:\S^+\to \N$ be the {\it prime-divisor function} defined by the formula
\[
\beta(S_1^{\Box\alpha_1}\Box S_2^{\Box\alpha_2} \Box\cdots \Box S_m^{\Box\alpha_m})=\alpha_1\alpha_2\cdots\alpha_m
\]
for all $S_1,\ldots,S_m\in \S^\Box$ and non-negative integers $\alpha_1,\ldots,\alpha_m$. Then $\beta(S)\le d(S)$ for each $S\in \S^+$, and thus $\beta$ has asymptotic mean value $1$ and asymptotic variance $0$.  
\end{ex}

\begin{ex}
Let $(\S,+,\Box,\partial)$ be a strictly monotonic arithmetical semiring for which $\cG_1^+$ holds. Let $\sigma_*:\S^+\to \R_{\ge 0}$ be the {\it divisor-sum function} defined by
\[
\sigma_*(S)=\sum_{S=S'\Box S''} \partial(S')
\]
for all $S\in \S^+$. Since 
\[
\sigma_*(S\Box S')=\sigma_*(S)(1+\partial(S'))+(1+\partial(S))\sigma_*(S')\le \sigma_*(S)\sigma_*(S')
\]
for all $\Box$-coprime $S,S'\in \S^+$ and
\[
\sigma_*(S^{\Box\alpha})=1+\partial(S)+\cdots+(\partial(S))^\alpha\le (1+\partial(S))^\alpha= (\sigma_*(S))^\alpha
\]
for all $S\in \S^\Box$ and $\alpha\in \N$, then $\sigma_*$ is submultiplicative. Moreover, Remark implies \ref{rem:monotonic} $n=O\left((\S^+\left(\left\lfloor\frac{n}{2}\right\rfloor\right)\right)$. Using the function $d$ introduced in Example \ref{ex:divisor}, we arrive at 
\[
(\sigma_*)_{\rm max}^+(n)\le (n+1) d^+_{\rm max}(n) = O(\S^+(n))\,.
\]
Since $\sigma_*$ has constant value $n+1$ on $\uS^\Box(n)$, Theorem \ref{thm:2} implies
\[
\sigma_*^+(n)=(n+1)\S^\Box(n)+o(\S^+(n))
\]
and thus $\sigma_*^+(n)\sim (n+1)\S^+(n)$.
\end{ex}

\begin{ex}
Let $(\S,+,\Box,\partial)$ be a strictly monotonic arithmetical semiring for which axiom $\cG_1^+$ holds. Let $\phi_*:\S^+\to \R_{\ge 0}$ be the {\it Euler-type function} that to each $S\in \uS^+(n)$ assigns the number $\phi_*(S)$ of elements of $\uS^+(n)$ that are $\Box$-coprime to $S$. It is easy to see that $\phi_*$ is submultiplicative and that $\phi_*$ restricts to the constant function $\S^+(n)-1$ on $\uS^\Box(n)$. Since $\S^+(n)-1$ is also the maximum of $\phi_*$ on $\uS^+(n)$, we conclude from Theorem \ref{thm:2} that 
\[
\frac{\phi_*(n)}{\S^+(n)}\sim \S^+(n)-1\sim \S^+(n)\,.
\]
\end{ex}

\end{document}